\documentclass[reqno]{amsart}

\usepackage{hyperref}
\usepackage{graphicx}

\newcommand{\pa}[2]{\frac{\partial #1}{\partial #2}}
\newcommand{\mr}[1]{\mathrm{#1}}
\newcommand{\mc}[1]{\mathcal{#1}}
\newcommand{\mf}[1]{\mathfrak{#1}}
\newcommand{\mb}[1]{\mathbf{#1}}
\newcommand{\mbb}[1]{\mathbb{#1}}

\newcommand{\xr}{\xrightarrow}

\DeclareMathOperator{\rank}{rank}

\DeclareMathOperator{\Ima}{Im}
\DeclareMathOperator{\Ker}{Ker}
\DeclareMathOperator{\Coker}{Coker}
\DeclareMathOperator{\Aut}{Aut}

\DeclareMathOperator{\diag}{diag}
\DeclareMathOperator{\lk}{lk}
\DeclareMathOperator{\St}{St}

\DeclareMathOperator{\Ad}{Ad}

\theoremstyle{plain}
\newtheorem{thm}{Theorem}[section]
\newtheorem{lem}[thm]{Lemma}

\newtheorem{prp}[thm]{Proposition}

\theoremstyle{definition}
\newtheorem{dfn}{Definition}[section]

\theoremstyle{remark}
\newtheorem{rmk}{Remark}[section]

\begin{document}

\title[Euclidean simplices and invariants of three-manifolds]{Euclidean simplices and invariants of three-manifolds:
a modification of the invariant for lens spaces}

\author{Evgeniy V. Martyushev}

\address{South Ural State University, 76 Lenin avenue, 454080 Chelyabinsk, Russia}

\email{mev@susu.ac.ru}

\thanks{This work has been partially supported by Russian Foundation for Basic Research, Grant no. 02-01-06003}

\begin{abstract}

We propose a three-manifold invariant based on the use of Euclidean metric values ascribed to the elements of manifold triangulation. We thus obtain
a nontrivial invariant that can, in particular, distinguish non-homeomorphic lens spaces.

\end{abstract}

\maketitle

\section*{Introduction}

It is well known that any two triangulations of a piecewise-linear three-manifold can be transformed into each other using the Pachner
moves~\cite{Pach}. If we construct an algebraic expression depending on some values ascribed to a manifold triangulation and invariant under these
moves, we get a value that is independent of the choice of triangulation. It is natural way to construct a PL-manifold invariant. For example, Euler
characteristic and Turaev-Viro type invariants can be proved to be invariant using exactly this method.

In this paper we present a new invariant of 3-dimensional PL-manifold based on Pachner's theorem. Its construction is naturally divided into three
main parts. First, on a given representation of the fundamental group we build a covering of a 3-manifold corresponding to the kernel of this
representation. Then, we map the covering space into 3-dimensional Euclidean space. In the last, algebraic part, we build an acyclic complex of
Euclidean geometric origin. The invariant of all the Pachner moves is expressed in terms of the torsion of that complex.

Such invariant was first constructed by I.G.~Korepanov in papers~\cite{Kor01, KM02} and calculated for some 3-manifolds, in particular, for a lens
space $L(p, q)$. Recall that at the end of paper~\cite{KM02} we proposed a general formula for our invariant for $L(p, q)$~\cite[formula
(4.1)]{KM02}. This formula gives, essentially, the squared Reidemeister torsion of $L(p, q)$. The main result of this paper is a proof of this
formula (theorem~\ref{thm:ILpq}).

The paper is organized as follows. In section~\ref{sec:adm} we prove an auxiliary theorem that will be used to prove the acyclicity of our complex
(theorem~\ref{thm:Gamma}). In section~\ref{sec:inv} we define an acyclic complex (theorem~\ref{thm:manacycl}) and prove that its torsion, divided by
some value, is a topological invariant, i.e., it is independent of the Pachner moves and some details of its construction (theorem~\ref{thm:inv}). In
section~\ref{sec:calc} we calculate our invariant for lens space $L(p, q)$ (theorem~\ref{thm:ILpq}) and show that for this case our invariant is
related to the Reidemeister torsion (remark~\ref{rmk:reid}). Finally, in section~\ref{sec:modif} we define and calculate a modified version of our
invariant for lens spaces (theorem~\ref{thm:ILpq2}). The modification uses a nontrivial action of the fundamental group on the acyclic complex.

\textbf{Acknowledgements.} I am grateful to I.G.~Korepanov for proposing me this problem and numerous helpful discussions.

\section{Admissible colorings and the map $\Gamma \colon \tilde{T} \to \mbb{R}^3$}\label{sec:adm}

Let $M$ be a connected closed orientable 3-dimensional manifold. Let $T$ be its \textit{triangulation}, i.e. a set $\Delta$ of pairwise disjoint
tetrahedra, together with a family of homeomorphisms, $\Phi$, mapping the set of tetrahedra faces to itself, so that the identification space
$\Delta/\Phi$ is homeomorphic to $M$. Note that this definition allows us to consider non-combinatorial manifold triangulations, that is, a 3-simplex
in it may not be determined uniquely by the set of its vertices. Moreover, simplex of any dimension is allowed to enter several times into the
boundary of higher dimensional simplex.

As is well known, for each subgroup $H$ of the fundamental group $\pi_1(M)$, there exists such a covering $p \colon \tilde{M} \to M$ that the induced
homomorphism $p_* \colon \pi_1(\tilde{M}) \to H$ is isomorphism. A covering corresponding to the trivial subgroup is called \textit{universal}. The
universal covering space $\tilde{M}$ is a \textit{simply connected} orientable 3-manifold. The fundamental group $\pi_1(M)$ acts on $\tilde{M}$
freely and transitively and $\tilde{M} / \pi_1(M) \cong M$.

Simplicial structure $T$ of manifold $M$ naturally induces a simplicial structure $\tilde{T}$ of $\tilde{M}$: every simplex from $T$ can be lifted to
$\tilde{T}$ (in general, in many different ways according to the action of $\pi_1(M)$).

Let us ascribe a real number $\lambda_{ij} > 0$ to every edge $e_{ij} \in T$, connecting vertices $i$ and $j$, so that the following inequality holds
for every tetrahedron of complex $T$ with edges $(e_{ij}, e_{ik}, e_{il}, e_{jk}, e_{jl}, e_{kl})$:
\begin{equation}\label{eq:CMdet}
\begin{vmatrix}
0 & 1 & 1 & 1 & 1 \\
1 & 0 & \lambda_{ij}^2 & \lambda_{ik}^2 & \lambda_{il}^2 \\
1 & \lambda_{ij}^2 & 0 & \lambda_{jk}^2 & \lambda_{jl}^2 \\
1 & \lambda_{ik}^2 & \lambda_{jk}^2 & 0 & \lambda_{kl}^2 \\
1 & \lambda_{il}^2 & \lambda_{jl}^2 & \lambda_{kl}^2 & 0
\end{vmatrix} > 0.
\end{equation}
Recall that this condition guarantees an existence of non-degenerate Euclidean tetrahedron with edge lengths $(\lambda_{ij}, \lambda_{ik},
\lambda_{il}, \lambda_{jk}, \lambda_{jl}, \lambda_{kl})$ and this tetrahedron is unique up to isometries, cf.~\cite{BerzheI}.

Further, let us ascribe a sign $+$ or $-$ to every tetrahedron $\tau \in T$. For each edge $a$, belonging to $i$th tetrahedron of complex $T$, we
denote by $\varphi_a^{(i)}$ its inner dihedral angle at this edge, taken with the sign ascribed to $i$th tetrahedron. This angle is a function of
$\lambda_{\ldots}$.

\begin{dfn}\label{dfn:deficit}
\textit{Defect angle} $\omega_a$ at an edge $a$ is a quantity
\[
\omega_a = - \sum_i \varphi_a^{(i)} \mod 2\pi,
\]
where the sum is taken over all tetrahedra sharing the edge $a$.
\end{dfn}

Let us say that we are given an \textit{admissible coloring} on $T$, if to every edge $e_{ij} \in T$ a number $\lambda_{ij}$ is ascribed so that the
condition~\eqref{eq:CMdet} holds, and to every tetrahedron from $T$ a sign is ascribed so that the defect angle at every edge vanishes.

Note that an admissible coloring on $T$ induces an admissible coloring on $\tilde{T}$ in a natural way: to every edge from the preimage of some
$e_{ij} \in T$ we ascribe the same number $\lambda_{ij}$, and to every tetrahedron from the preimage of some $\tau \in T$ we ascribe the same sign as
to $\tau$.

\begin{thm}\label{thm:Gamma}
Let $\tilde{T}$ be a simplicial complex for the universal covering space $\tilde{M}$. Let an admissible coloring be given on $\tilde{T}$. Then, there
exists a continuous map $\Gamma \colon \tilde{T} \to \mbb{R}^3$ such that $l_{ij} = \lambda_{ij}$, $\forall i, j$, where $l_{ij}$ is the length of
$\Gamma(e_{ij})$. Any other map $\Gamma' \colon \tilde{T} \to \mbb{R}^3$ for the same admissible coloring can be obtained from $\Gamma$ by an
orientation preserving isometry of $\mbb{R}^3$.
\end{thm}

\begin{proof}
Let $\tau_1$ be a certain tetrahedron of $\tilde{T}$ and let $(\lambda_{12}, \lambda_{13}, \lambda_{14}, \lambda_{23}, \lambda_{24}, \lambda_{34})$
be an admissible coloring of its edges. Denote by $\Gamma_1 \colon \tau_1 \to \mbb{R}^3$ an arbitrary embedding of this tetrahedron in $\mbb{R}^3$
such that $l_{ij} = \lambda_{ij}$. Let $\tau_2$ be an adjacent tetrahedron, i.e. $\tau_1 \cap \tau_2$ is 2-dimensional simplex. Then, there exists an
embedding $\Gamma_2 \colon \tau_2 \to \mbb{R}^3$ such that $\Gamma_2 = \Gamma_1$ on $\tau_2 \cap \tau_1$. This embedding is unique: the tetrahedra
$\Gamma_1(\tau_1)$ and $\Gamma_2(\tau_2)$ are on the same side or on the different sides from their common 2-face, if the signs of $\tau_1$ and
$\tau_2$ are the same or different respectively. We say that $\Gamma_2$ continue $\Gamma_1$ by means of $\tau_2$.

Let us take an arbitrary tetrahedron $\tau_i \in \tilde{T}$ and consider a sequence of pairwise adjacent tetrahedra starting from $\tau_1$, that is,
a sequence $s_i = (\tau_1, \ldots, \tau_i)$, where $\tau_j \cap \tau_{j+1}$ is a 2-simplex for each consecutive pair $(\tau_j, \tau_{j+1})$ and for
all $j \in \{1, \ldots, i\}$. Finally, define $\Gamma_{s_i} \colon \tau_i \to \mbb{R}^3$ as a continuation of $\Gamma_1$ by means of the sequence
$s_i$.

\begin{lem}\label{lem:indep}
The map $\Gamma_{s_i}$ is independent of a sequence $s_i$, i.e., if $s'_i$ is any other sequence of adjacent tetrahedra joining $\tau_1$ and
$\tau_i$, then $\Gamma_{s_i} = \Gamma_{s'_i} = \Gamma_i$.
\end{lem}

\begin{proof}
We have a sequence $s_i$ connecting $\tau_1$ and $\tau_i$. It can be geometrically imagined as a broken line whose every straight segment connects
the barycenters of two adjacent tetrahedra. We must show that if $s'_i$ is any other sequence connecting $\tau_1$ and $\tau_i$, then $\Gamma_1 =
\Gamma'_1$, where $\Gamma'_1 \colon \tau_1 \to \mbb{R}_3$ is a continuation of $\Gamma_1$ by means of a closed path $s_i \circ (s'_i)^{-1}$.

Recall that we are given an admissible coloring on $\tilde{T}$. It follows that all the defect angles are zeros. Let $n$ tetrahedra $\tau_{k_1},
\ldots, \tau_{k_n}$ share an only edge $e$, where $\tau_{k_i}$ and $\tau_{k_{i+1}}$ share a face for all $i$. Denote by $\alpha$ a closed path
$(\tau_{k_1}, \ldots, \tau_{k_n}, \tau_{k_1})$. Since the defect angle at $e$ is zero, then $\Gamma_{k_1} = \Gamma'_{k_1}$, where $\Gamma'_{k_1}
\colon \tau_{k_1} \to \mbb{R}_3$ is a continuation of $\Gamma_{k_1}$ by means of $\alpha$.

Using Seifert-Van Kampen theorem, one can show that a closed path $s_i \circ (s'_i)^{-1}$ may be represented as a finite product of paths like
$\alpha$. Hence lemma~\ref{lem:indep} follows.
\end{proof}

Thus, for each tetrahedron $\tau_i \in \tilde{T}$ we have an embedding $\Gamma_i \colon \tau_i \to \mbb{R}^3$ so that two embeddings $\Gamma_i$ and
$\Gamma_j$ agree on $\tau_i \cap \tau_j, \forall i, j$. Now we define $\Gamma \colon \tilde{T} \to \mbb{R}^3$ by
\[
\Gamma(x) = \Gamma_i(x) \Longleftrightarrow x \in \tau_i, \forall i.
\]
In order to complete the proof of theorem~\ref{thm:Gamma}, we need the following well-known result from general topology.
\begin{lem}\label{lem:pasting}
Let $X = A \cup B$, where $A$ and $B$ are the closed subspaces of $X$. Let $f \colon A \to Y$ and $g \colon B \to Y$ be continuous. If $f = g$ on $A
\cap B$, then a function $h$ defined by
\[
h(x) =
\begin{cases}
f(x), & x \in A, \\
g(x), & x \in B,
\end{cases}
\]
is continuous.
\end{lem}

Lemma~\ref{lem:pasting} and definition of $\Gamma$ imply that it is continuous.

The last statement of the theorem concerns the arbitrariness in the choice of embedding $\Gamma_1$. Let $\Gamma'_1$ be another embedding of $\tau_1$
in $\mbb{R}^3$. Denote by $E(3)$ the group of orientation preserving motions of 3-dimensional Euclidean space $\mbb{R}^3$. Then, there is an element
$h \in E(3)$ such that $\Gamma'_1(\tau_1) = h\Gamma_1(\tau_1)$. It follows from the construction of $\Gamma$ that $\Gamma'(\tilde{T}) =
h\Gamma(\tilde{T})$. Hence, $\Gamma$ is independent of a choice of the first tetrahedron $\tau_1$ up to multiplying by an element of $E(3)$. This
completes the proof of theorem~\ref{thm:Gamma}.
\end{proof}

Let $\tilde{\tau}$ and $\tilde{\tau}'$ be preimages of an only tetrahedron $\tau \in T$ with respect to the covering projection. Then, there is an
element $g \in \pi_1(M)$ such that $\tilde{\tau}' = g\tilde{\tau}$. Since the edge lengths of $\tilde{\tau}$ and $\tilde{\tau}'$ are the same, there
exists an element $h \in E(3)$ such that $\Gamma(\tilde{\tau}') = h\Gamma(\tilde{\tau})$. Thus we have a mapping $\rho \colon g \mapsto h$.
Obviously, $\rho(g_1 g_2) = \rho(g_1) \rho(g_2)$, hence $\rho \colon \pi_1(M) \to E(3)$ is a homomorphism.

Denote $\Gamma' = h \circ \Gamma$, where $h \in E(3)$. Then, $\rho' = h^{-1} \circ \rho \circ h$, where $\rho'$ is a representation corresponding to
$\Gamma'$. Recall that the representations $\rho, \rho' \colon \pi_1 (M) \to E(3)$ are called \textit{equivalent} in this case. So, the following
proposition holds.
\begin{prp}\label{prp:inject}
For any admissible coloring of $\tilde{T}$, there is a unique class $[\rho]$ of equivalent representations $\pi_1(M)$ in $E(3)$.
\end{prp}

Conversely, given a representation
\begin{equation}\label{eq:repr}
\rho \colon \pi_1(M) \to E(3),
\end{equation}
one can construct a continuous map $\Gamma_{\rho} \colon \tilde{T} \to \mbb{R}^3$ as follows.

First of all, note that all the tetrahedra in the triangulation of orientable space $\tilde{M}$ can be oriented \textit{consistently}, i.e., for
every tetrahedron, we can order its vertices up to even permutations. For example, we suppose the vertex orderings $ABCD$ and $EABC$ of two adjacent
tetrahedra to be consistent.

\begin{dfn}\label{dfn:family}
\textit{A fundamental family of simplices} in $\tilde{T}$ is such a family $\mc{F}$ of simplices of $\tilde{T}$ that over each simplex of $T$ lies
exactly one simplex of this family.
\end{dfn}

Let us fix a consistent orientation on the set of all tetrahedra from $\tilde{T}$ and also the fundamental family $\mc{F}$ for $\tilde{T}$.

Let $v$ be a vertex belonging to $T$ and let $\tilde{v} \in \mc{F}$ be its preimage in the fundamental family. Then we place the orbit of
$\tilde{v}$, i.e. the set $\{g \tilde{v} \mid g \in \pi_1(M) \}$, in $\mbb{R}^3$ so that if $\tilde{v}_2 = g \tilde{v}_1$, where $g \in \pi_1(M)$,
then $\Gamma_{\rho}(\tilde{v}_2) = \rho (g) \Gamma_{\rho}(\tilde{v}_1)$. Thus, to each vertex of $T$ we assign a set of points in $\mbb{R}^3$ whose
coordinates are related by the motions of $\mbb{R}^3$.

\begin{rmk}\label{rmk:genpos}
The vertices from $\mc{F}$ are mapped in $\mbb{R}^3$ in arbitrary way. However, we require the configuration of these vertices in $\mbb{R}^3$ to obey
the following conditions of general position:
\begin{itemize}
\item
the volumes of all tetrahedra from $\Gamma(\mc{F})$ are nonzero;
\item
the rank of $\left(\pa{\omega_i}{l_j}\right)$ possesses the maximal value at the point $\omega_i=0$, $\forall i$; here, $\omega_i$ is the defect
angle at $i$th edge, $l_j$ is the length of $j$th edge, $i$ and $j$ run over all the edges from $\mc{F}$.
\end{itemize}
\end{rmk}

Under $\Gamma_{\rho}$, each simplex $\sigma^k \in \tilde{T}$ of nonzero dimension $k$ is mapped into the convex linear shell of its vertex images in
$\mbb{R}^3$. That is, if $v_0$, \ldots, $v_k$ are the vertices of $\sigma^k$, then $\Gamma_{\rho}(\sigma^k) = \sum\limits_{i = 0}^k m_i
\Gamma_{\rho}(v_i)$, where $m_0 + \ldots + m_k = 1$ and $m_i \geq 0, \, \forall i$.

Hence, every edge $e \in \tilde{T}$ is assigned the length of $\Gamma_{\rho}(e) \in \mbb{R}^3$, and every 3-simplex $\tau \in \tilde{T}$ is assigned
the oriented volume of $\Gamma_{\rho}(\tau)$. Similarly, we can also speak about dihedral angles between 2-simplices in the triangulation.

Under $\Gamma_{\rho}$, the orientation of a tetrahedron either coincides with the orientation of $\mbb{R}^3$ or does not. In the first case we take
the volume of the tetrahedron and all its dihedral angles with the sign $+$, in the second case with the sign $-$.

Thus, given representation $\rho$, we have found an admissible coloring on $\tilde{T}$: for every edge we ascribe the length of its image in
$\mbb{R}^3$, for every tetrahedron we ascribe the sign of its volume in $\mbb{R}^3$.

Before constructing the invariant, we need to recall some formulas of Euclidean geometry.

\subsection{Formulas of Euclidean geometry}

We are going to write out some formulas for derivatives $\pa{\omega_a}{l_b}$, taken under condition that the lengths of all edges, except $b$, are
fixed. These formulas can be easily proved by direct calculation using definitions and elementary geometric relations like cosine rule,
see~\cite{KM02} for details.

Clearly, we obtain nonzero derivatives only if $a$ and $b$ lie in a single tetrahedron. Consider several cases.
\begin{description}
\item[1st case]
Edges $a = AB$ and $b = CD$ are skew in a tetrahedron $ABCD$ and there is no other tetrahedron in the triangulation containing both edges $a$ and
$b$:
\begin{equation}\label{eq:skresch}
\pa{\omega_{AB}}{l_{CD}} = - \frac{l_{AB} l_{CD}}{V_{ABCD}},
\end{equation}
where $V_{ABCD}$ is the volume of $ABCD$ multiplied by 6.

\item[2nd case]
Edges $a = AB$ è $b = BC$ belong to the common face of two adjacent tetrahedra and, again, there is no other tetrahedron in the triangulation
containing both edges $a$ and $b$:
\begin{equation}\label{eq:peresek}
\pa{\omega_{AB}}{l_{BC}} = l_{AB} l_{BC} \, \frac {V_{CAED}}{V_{ABCD} V_{EABC}}.
\end{equation}

\item[3rd case]
Edge $a = b = DE$ is common for exactly three tetrahedra $ABED$, $BCED$ and $CAED$, and no one of these tetrahedra contain this edge more than one
time:
\begin{equation}\label{eq:sovpad}
\pa{\omega_{DE}}{l_{DE}} = - l_{DE}^2 \, \frac{V_{ABCD} V_{EABC}}{V_{ABED} V_{BCED} V_{CAED}}.
\end{equation}

\item[4th case]
Edge $a = b = DE$ is common for $>3$ tetrahedra, and again no one of these tetrahedra contain this edge more than one time.

We write out this formula only for the case of four tetrahedra (the generalization for greater number of tetrahedra is obvious):
\begin{equation}\label{eq:4tet}
\pa{\omega_{DE}}{l_{DE}} = - \frac {l_{DE}^2}{6} \left(\frac{V_{ABCD} V_{EABC}}{V_{ABED} V_{BCED} V_{CAED}} + \frac{V_{AFBD} V_{EAFB}}{V_{AFED}
V_{FBED} V_{BAED}}\right).
\end{equation}
Here, the right-hand side consists of two terms: the first one is the same as in~\eqref{eq:sovpad}, and the second one is obtained by replacing $B
\to F$ and $C \to B$.

\item[5th case]
Actually, our definition of triangulation admits such gluings between faces of tetrahedra where different combinations of above mentioned cases
appear. Namely, the edge $b$ may
\begin{itemize}
\item
lie opposite to edge $a$ (1st case) in more than one tetrahedron;
\item
belong to the same 2-dimensional face as $a$ (2nd case), and again there may be more than one such faces;
\item
coincide with $a$ (as in the 3rd and 4th case),
\end{itemize}
and these possibilities (as examples show) do not exclude one another. Thus, there may be several ``ways of influence'' of $dl_b$ on $d\omega_a$.
Clearly, the terms corresponding to these ways sum together.
\end{description}

\section{Acyclic complex and invariant of 3-manifold}\label{sec:inv}

\subsection{Generalities on acyclic complexes and their torsions}

We shortly remind basic definitions from the theory of algebraic complexes, see~\cite{Tur01} for details.

Let $C_0$, $C_1$, \ldots , $C_n$ be finite-dimensional $\mbb{R}$-vector spaces. We suppose that each $C_i$ is based, that is has distinguished basis.
Then, linear mapping $f_i \colon C_{i + 1} \to C_i$ can be identified with matrix.

\begin{dfn}\label{dfn:cmplx}
The sequence of vector spaces and linear mappings
\begin{equation}\label{eq:cmplx}
C = (0 \xr{} C_n \xr{f_{n-1}} C_{n-1} \xr{} \ldots \xr{} C_1 \xr{f_0} C_0 \xr{} 0)
\end{equation}
is called a \textit{complex} if $\Ima f_i \subset \Ker f_{i-1}$ for all $i = 1, \ldots , n - 1$. This condition is equivalent to $f_{i-1} f_i = 0$
for all $i$.
\end{dfn}

\begin{dfn}\label{dfn:homol}
The space $H_i(C) = \Ker f_{i-1} / \Ima f_i$ is called the \textit{$i$th homology} of the complex $C$.
\end{dfn}

\begin{dfn}\label{dfn:acycl}
The complex $C$ is said to be \textit{acyclic} if $H_i(C) = 0$ for all $i$. This condition is equivalent to $\rank f_{i-1} = \dim C_i - \rank f_i$
for all $i$.
\end{dfn}

Suppose that the sequence~\eqref{eq:cmplx} is an acyclic complex. Let $\mc{C}_i$ be an ordered set of basis vectors in $C_i$ and let $\mc{B}_i
\subset \mc{C}_i$ be a subset of basis vectors belonging to the space $\Ima f_i$.

Denote by ${}_{\mc{B}_i} f_i$ a nondegenerate transition matrix from the basis in space $\Coker f_{i+1} = C_{i+1}/\Ima f_{i+1}$ to the basis in space
$\Ima f_i$. By acyclicity, such a matrix really exists. Hence, ${}_{\mc{B}_i} f_i$ is a principal minor of the matrix $f_i$ obtained by striking out
the rows corresponding to vectors of $\mc{B}_{i+1}$ and the columns corresponding to vectors of $\mc{C}_i \setminus \mc{B}_i$.

\begin{dfn}\label{dfn:tors}
A quantity
\begin{equation}\label{eq:torsion}
\tau(C) = \prod\limits_{i = 0}^{n-1} (\det {}_{\mc{B}_i} f_i)^{(- 1)^{i+1}}
\end{equation}
is called the \textit{torsion} of acyclic complex $C$.
\end{dfn}

\begin{thm}[\cite{Tur01}]
Up to a sign, $\tau(C)$ does not depend on the choice of subsets $\mc{B}_i$.
\end{thm}

\begin{rmk}\label{rmk:changing}
The torsion $\tau(C)$ does depend on the distinguished basis of $C_i$. If one performs change-of-basis transformation in every space $C_i$ with
nondegenerate matrix $A_i$, then the torsion $\tau(C)$ is multiplied by
\[
\prod\limits_{i = 0}^n (\det A_i)^{(-1)^{i+1}}.
\]
\end{rmk}

From now on we suppose the fundamental group $\pi_1(M)$ to be finite. Since any real representation of a finite group is equivalent to the orthogonal
one, then $\Ima \rho$ is actually a subgroup in $\mr{SO}(3)$.

\begin{dfn}
For a given representation $\rho \colon \pi_1(M) \to E(3)$, the \textit{centralizer} of the subgroup $\Ima \rho$ is a set of elements $g \in E(3)$ so
that $g \rho(h) g^{-1} = \rho(h)$ for all $h \in \pi_1(M)$. The centralizer is a subgroup in $E(3)$, and we will denote it by $E(3)_\rho$. We denote
by $\mf{e}(3)_\rho$ the Lie algebra of $E(3)_{\rho}$. So,
\[
\mf{e}(3)_\rho = \{u \in \mf{e}(3) \mid \Ad_{\rho(h)} u = u, \forall h \in \pi_1(M)\},
\]
where $\Ad_\rho = \Ad \circ \rho \colon \pi_1(M) \to \Aut(\mf{e}(3))$ and $\Ad \colon g \mapsto \Ad_g$ is the \textit{adjoint} representation.
\end{dfn}

Let us take a geometrically natural basis of the whole Lie algebra $\mf{e}(3)$ of $E(3)$. This basis consists of three translations $dx$, $dy$ and
$dz$ along the Cartesian axes and three rotations $d\varphi_x$, $d\varphi_y$ and $d\varphi_z$ around these axes. Then, we also have a basis in
$\mf{e}(3)_\rho$ induced by restriction.

\begin{dfn}
A representation $\rho \colon \pi_1 (M) \to E(3)$ is called \textit{abelian} (resp. \textit{trivial}) if its image is abelian (resp. trivial)
subgroup of $E(3)$.
\end{dfn}

\begin{prp}\label{prp:e3rho}
Let $\mf{e}(3)_2$ be a 2-dimensional subalgebra of $\mf{e}(3)$ generated by the differential of screw motion along a fixed axis in $\mbb{R}^3$ and
let $\mf{o}$ be the trivial algebra Lie. Then,
\begin{equation}
\mf{e}(3)_\rho =
\begin{cases}
\mf{e}(3), & \hbox{\text{if $\rho$ is trivial};} \\
\mf{e}(3)_2, & \hbox{\text{if $\rho$ is abelian};} \\
\mf{o}, & \hbox{\text{if $\rho$ is non abelian}.}
\end{cases}
\end{equation}
\end{prp}

\begin{proof}
If $\rho$ is trivial, then the centralizer of $\Ima \rho$ is $E(3)$, and hence, $\mf{e}(3)_\rho = \mf{e}(3)$.

Let $\rho$ be abelian. Then, $\Ima \rho$ consists of rotations around the only axis $z$ in $\mbb{R}^3$. There exist two types of motions in $E(3)$
which commute with every element of $\Ima \rho$: a rotation around $z$ and translation along this axis. Hence, $\mf{e}(3)_\rho$ is generated by the
differential of a screw motion along $z$.

Let $\rho$ be non abelian. In this case, $\Ima \rho$ consists of rotations around at least two nonparallel axes. The existence of such axes fixes the
coordinate system in $\mbb{R}^3$. Therefore, $E(3)_\rho = 1$ and $\mf{e}(3)_\rho = \mf{o}$.
\end{proof}

Denote by $N_i$ the full amount of $i$-dimensional simplices in $\mc{F}$.

We define a $3N_0$-dimensional vector space $\mbb{R}^{3N_0}_x$ as a space of all mappings of the vertices from $\mc{F}$ to $\mbb{R}^3$. Denote by
$(dx) = T_{\Gamma} \mbb{R}^{3N_0}_x$ the tangent space to $\mbb{R}^{3N_0}_x$ at $\Gamma$. The basis in $(dx)$ is generated by a choice of coordinate
system in $\mbb{R}^3$. For convenience, we choose the \textit{Cartesian} coordinates.

A linear mapping $f_1 \colon \mf{e}(3)_\rho \to (dx)$ is defined as follows. Under $f_1$, the six generators $dx$, $dy$, $dz$, $d\varphi_x$,
$d\varphi_y$ and $d\varphi_z$ of $\mf{e}(3)$ are mapped into the differentials of a vertex $A$ according to the evident formula:
\begin{equation}\label{eq:dxdydz}
\begin{pmatrix}
dx_A \\ dy_A \\ dz_A
\end{pmatrix} =
\begin{pmatrix}
0 & d\varphi_z & - d\varphi_y \\
- d\varphi_z & 0 & d\varphi_x \\
d\varphi_y & - d\varphi_x & 0
\end{pmatrix}
\begin{pmatrix}
x_A \\ y_A \\z_A
\end{pmatrix} +
\begin{pmatrix}
dx \\ dy \\ dz
\end{pmatrix},
\end{equation}
where $x_A$, $y_A$ and $z_A$ are the coordinates of $A$.

Let $e \in \mc{F}$ be any edge of the fundamental family. Let $L_e = \frac{1}{2} l_e^2$, where $l_e$ is the Euclidean length of $\Gamma(e) \in
\mbb{R}^3$. Define $N_1$-dimensional vector space $\mbb{R}^{N_1}_L$ as a space of all mappings of the edges from $\mc{F}$ to $\mbb{R}$. The subscript
$L$ stresses that the basis of this space consists of the vectors $(L_1,\ldots, L_{N_1})$. Denote by $(dL) = T_{\Gamma} \mbb{R}^{N_1}_L$ the tangent
space to $\mbb{R}^{N_1}_L$ at $\Gamma$.

Consider a linear mapping $f_2 \colon (dx) \to (dL)$. By definition, $f_2$ sends the given differentials of two vertices $A$ and $B$ to
\begin{equation}\label{eq:dL}
dL_{AB} = (x_B - x_A) (dx_B - dx_A) + (y_B - y_A) (dy_B - dy_A) + (z_B - z_A) (dz_B - dz_A).
\end{equation}

Since the bases of $\mf{e}(3)_\rho$, $(dx)$ and $(dL)$ are distinguished, then the mappings $f_1$ and $f_2$ are uniquely determined by their matrices
which we denote by the same letters.

Put $\Omega_e = \frac{\omega_e}{l_e}$, where $\omega_e$ is the defect angle at $\Gamma(e)$ (definition~\ref{dfn:deficit}). Similarly as we have
defined the space $(dL)$, let us define the vector space $(d\Omega) = T_{\Gamma} \mbb{R}^{N_1}_{\Omega}$: its basis consists of the vectors
$(d\Omega_1,\ldots, d\Omega_{N_1})$.

Let us number the edges of $\mc{F}$ in a certain order and consider the matrix $f_3 = \left(\pa{\Omega_a}{L_b}\right)$, where $a$ and $b$ run over
the edges of $\mc{F}$ and the partial derivatives are taken at $\omega_1 = \ldots = \omega_{N_1} = 0$. The matrix elements of $f_3$ are defined by
formulas~\eqref{eq:skresch}, \eqref{eq:peresek}, \eqref{eq:sovpad} and~\eqref{eq:4tet}. In particular, these formulas imply that $f_3$ is symmetric.
We can obtain the same result directly.
\begin{lem}
The matrix $f_3$ is symmetric, i.e. $f_3^T = f_3$.
\end{lem}

\begin{proof}
Since $f_3 = \left(\frac{1}{l_a l_b}\pa{\omega_a}{l_b}\right)$, it suffices to prove the symmetric form of matrix $\left(\pa{\omega_a}{l_b}\right)$.
Let us find the Hessian of the function $H(l_1,\ldots, l_{N_1}) = \sum\limits_k l_k \omega_k$ at the point $\omega_1 = \ldots = \omega_{N_1} = 0$:
\[
\pa{}{l_b}\left(\pa{H}{l_a}\right) = \pa{}{l_b}\left(\omega_a + \sum\limits_{k=1}^{N_1} l_k \pa{\omega_k}{l_a}\right) = \pa{\omega_a}{l_b}.
\]
The second equality holds since $\sum\limits_{k=1}^{N_1} l_k d \omega_k = 0$ (it is a consequence of well-known Schl\"{a}fli's differential
identity~\cite[p.~281--295]{Mil}). Hence, the claim follows.
\end{proof}

Denote by $\mf{e}(3)^*_\rho$ and $(dx)^*$ vector spaces whose distinguished bases are \textit{dual} to the ones of $\mf{e}(3)_\rho$ and $(dx)$
respectively. This means that if, say, $d = (d_1, \ldots, d_k)$ is a basis of $\mf{e}(3)_\rho$, then a basis $d^* = (d^*_1, \ldots, d^*_k)$ of
$\mf{e}(3)^*_\rho$ is defined by
\[
d^*_i(d_j) =
\begin{cases}
1, & i = j, \\
0, & i \neq j.
\end{cases}
\]

\begin{thm}\label{thm:manacycl}
The sequence of vector spaces and linear mappings
\begin{equation}\label{eq:seq}
0 \xr{} \mf{e}(3)_\rho \xr{f_1} (dx) \xr{f_2} (dL) \xr{f_3 = f_3^T} (d\Omega) \xr{- f_2^T} (dx)^* \xr{f_1^T} \mf{e}(3)^*_\rho \xr{} 0
\end{equation}
is an acyclic complex (definition~\ref{dfn:acycl}).
\end{thm}

\begin{rmk}\label{rmk:acycl}
Recall that we restrict ourselves to the case of finite fundamental group. An example of space $S^2 \times S^1$ shows that for manifolds with an
infinite fundamental group there can exist a representation $\rho$ such that~\eqref{eq:seq} is not acyclic. However, it turns out that for this case
one can define an acyclic complex like~\eqref{eq:seq} as well (see~\cite{KII}).
\end{rmk}

\noindent\textit{Proof of theorem~\ref{thm:manacycl}.}

\noindent\textbf{Acyclicity at $\mf{e}(3)_\rho$} (= injectivity of $f_1$) holds because the vertices of the fundamental family are in general
position (remark~\ref{rmk:genpos}).

\bigskip

\noindent\textbf{Acyclicity at $(dx)$}. The elements of $\mf{e}(3)_\rho$ generate global infinitesimal translations of the vertices from $\mc{F}$
which do not change the lengths, according to~\eqref{eq:dxdydz}. It follows that $\Ima f_1 \subset \Ker f_2$. Conversely, if the vertex coordinates
take arbitrary infinitesimal increments so that the differentials of all the lengths remain zero, then this corresponds to an element of the algebra
$\mf{e}(3)$. It can be shown by direct calculation that if $h$ is such an element of $\pi_1(M)$ that $u \neq \Ad_{\rho(h)} u$, then any nonzero
element $u \in \mf{e}(3) / \mf{e}(3)_\rho$ generates infinitesimal increments of the vertices for which $dL_{A, hA} \neq 0$, where $A$ is a vertex
from $\mc{F}$. Clearly, it is impossible if $dL_e = 0$ for every $e$ in $\mc{F}$. Hence, $\Ker f_2 \subset \Ima f_1$ and the acyclicity at $(dx)$
follows.

\bigskip

\noindent\textbf{Acyclicity at $(dL)$}. The elements of factor-space $(dx)/ \mf{e}(3)_\rho$ generate infinitesimal increments of the edge lengths,
but all the defect angles remain zero, i.e. $\Ima f_2 \subset \Ker f_3$. To prove the converse, we must show that if all $L_j$ take infinitesimal
increments $dL_j$ so that all $d\Omega_i = 0$, then $dL_j$'s are generated by some elements of $(dx)$. Let us first prove the ``finite'' version of
this statement.

\begin{lem}\label{lem:finite}
Let the lengths of the edges from $\mc{F}$ take small but finite increments $\Delta l_j$ so that all $\Delta\omega_i = 0$. Then, there are such
points in $\mbb{R}^3$ with coordinates $(x'_k, y'_k, z'_k)$ that
\begin{equation}\label{eq:DeltaL}
\Delta L_j = l_j \Delta l_j = (x_k - x_s) (\Delta x_k - \Delta x_s) + (y_k - y_s) (\Delta y_k - \Delta y_s) + (z_k - z_s) (\Delta z_k - \Delta z_s),
\end{equation}
where $\Delta x_k = x'_k - x_k$.
\end{lem}

\begin{proof}
If $\Delta l_j$ are small enough, then all the determinants from~\eqref{eq:CMdet} remain positive. With a condition $\Delta \omega_i = 0, \forall i$
this means that $l_j + \Delta l_j$ give an admissible coloring on $T$ and induced admissible coloring on $\tilde{T}$. By theorem~\ref{thm:Gamma},
there exists a continuous map $\Gamma' \colon \tilde{T} \to \mbb{R}^3$ such that $l'_j = l_j + \Delta l_j$, $\forall j$, where $l'_j$ is the length
of $\Gamma'(e_j)$ of some edge $e_j$. Then, $\Gamma'$ maps all the vertices from the fundamental family to points in $\mbb{R}^3$ which are the ones
we looked for.
\end{proof}

This lemma implies that for the sequence
\begin{equation}
\mbb{R}^{3N_0}_x \xr{F_2} \mbb{R}^{N_1}_L \xr{F_3} \mbb{R}^{N_1}_\Omega,
\end{equation}
the inclusion $\Ker F_3 \subset \Ima F_2$ holds. Note that $f_2 = dF_2$ and $f_3 = dF_3$. Therefore, taking the differential of this sequence and
taking into account the conditions of general position from remark~\ref{rmk:genpos}, we get $\Ker f_3 \subset \Ima f_2$, hence the acyclicity at
$(dL)$ follows.

\bigskip

\noindent\textbf{Acyclicity at other spaces}. Recall that $\Ker f_3 = \Ima f_2$ is equivalent to
\[
\rank f_3 = \dim (dL) - \rank f_2.
\]
Taking into account that $\dim (dL) = \dim (d\Omega)$ and $\rank (-f_2^T) = \rank f_2$, we get
\[
\rank (-f_2^T) = \dim (d\Omega) - \rank f_3,
\]
that means the acyclicity at $(d\Omega)$. Likewise, the acyclicity at $(dx)^*$ and $\mf{e}(3)_\rho^*$ follows from the acyclicity at $(dx)$ and
$\mf{e}(3)_\rho$ respectively. Theorem~\ref{thm:manacycl} is proven. \qed \medskip

Denote by $\mc{C}_0$, $\mc{C}_1$ and $\mc{C}_2$ arbitrary ordered sets of basis vectors in the spaces $\mf{e}(3)_\rho$, $(dx)$ and $(dL)$
respectively. Let $\mc{B}_i \subset \mc{C}_i$ be a subset of basis vectors belonging to $\Ima f_i$. Denote by ${}_{\mc{B}_i} f_i$ such a principal
minor of the matrix $f_i$ that its rows correspond to the vectors from $\mc{C}_{i-1} \setminus \mc{B}_{i-1}$ and its columns correspond to the
vectors from $\mc{B}_i$. Let $\mc{B}_3 = \mc{C}_2 \setminus \mc{B}_2$. Then, according to definition~\ref{dfn:tors}, the torsion of the
complex~\eqref{eq:seq} looks like
\[
\tau = \frac{(\det {}_{\mc{B}_2} f_2)^2 \, (-1)^{\rank f_2}}{\det {}_{\mc{B}_3} f_3 \, (\det {}_{\mc{B}_1} f_1)^2}.
\]
By the acyclicity,
\[
(-1)^{\rank f_2} = (-1)^{\dim (dx) - \dim \mf{e}(3)_\rho} = (-1)^{3N_0} = (-1)^{N_0},
\]
where the second equality is due to proposition~\ref{prp:e3rho}. It follows that
\begin{equation}\label{eq:tors}
\tau = \frac{(\det {}_{\mc{B}_2} f_2)^2 \, (-1)^{N_0}}{\det {}_{\mc{B}_3} f_3 \, (\det {}_{\mc{B}_1} f_1)^2}.
\end{equation}
Besides, the form of~\eqref{eq:seq} implies that the torsion $\tau$ is independent of the ordering on $\mc{C}_i, \forall i$.

Put
\begin{equation}\label{eq:inv}
I_\rho(M) =  \frac{\tau}{\prod (-V)},
\end{equation}
where $V$ is the volume of tetrahedron multiplied by 6, and the product is taken over all tetrahedra of $\mc{F}$.

\begin{thm}\label{thm:inv}
For a given representation $\rho$, the quantity $I_\rho(M)$ is an invariant of the manifold $M$.
\end{thm}

\noindent\textit{Proof} is by a series of lemmas. We have to show that $I_\rho(M)$ is independent of various details of its construction. First, it
does not depend on the choice of the fundamental family $\mc{F}$. Second, using formula~\eqref{eq:sovpad} and Pachner's theorem, one can show that
$I_\rho(M)$ does not depend on the choice of the triangulation of $M$ (see~\cite{Kor01} for details of proof). Finally, $I_\rho(M)$ is independent of
$\Gamma$, i.e. it does not matter in which points of Euclidean space we place the vertices from the fundamental family (see remark~\ref{rmk:genpos}).

\begin{lem}\label{lem:fund}
$I_\rho(M)$ does not depend on the choice of the fundamental family $\mc{F}$.
\end{lem}

\begin{proof}
First, note that the corresponding simplices from the different fundamental families can be obtained from each other by an action of a certain
element of $\pi_1(M)$. It follows that the lengths, dihedral angles and tetrahedron volumes do not change under replacing the fundamental family with
another one. Therefore, the bases in $(dx)$ and $(dx)^*$ can only change and they must stay dual to each other.

Denote by $A, B \in \mr{SO}(3N_0)$ the change-of-basis matrices in the spaces $(dx)$ and $(dx)^*$ respectively. Then, $\widetilde{f_2} = f_2 A$ and
$\widetilde{- f_2^T} = - B^{-1} f_2^T = - \widetilde{f_2}^T = - A^T f_2^T$. Hence, $B = (A^T)^{-1}$. That is, the torsion of~\eqref{eq:seq} is
divided by $(\det A)^2 = 1$ (see~remark~\ref{rmk:changing}) and the claim follows.
\end{proof}

\begin{lem}[\cite{Kor01}]\label{lem:pach}
$I_\rho(M)$ does not depend on the triangulation of $M$.
\end{lem}

\begin{lem}\label{lem:proizvol}
$I_\rho(M)$ is independent of $\Gamma$, i.e., it does not matter in which points of Euclidean space we place the vertices from $\mc{F}$.
\end{lem}

\begin{proof}
Let $v$ be an arbitrary vertex of $T$ and $\tilde{v} \in \mc{F}$. We are going to show that using the Pachner moves one can move $O =
\Gamma(\tilde{v})$ in any other point of $\mbb{R}^3$ in such a way that all the remaining vertices from $\Gamma(\mc{F})$ do not change their
positions. Then, the claim will follow from lemma~\ref{lem:pach}.

Recall that a union of all simplices containing $O$ as a vertex is called the \textit{star} and is denoted by $\St(O)$. A union of all simplices of
$\St(O)$ which do not contain $O$ is called the \textit{link} of the vertex $O$ and is denoted by $\lk(O)$.

First of all, let us make a \textit{barycentric} subdivision (which can be represented as a finite sequence of the Pachner moves) of the initial
triangulation $T$ so as to reach a \textit{combinatorial} triangulation where the link of every vertex is a triangulated 2-sphere.

Further, let vertices $A, B, C$ and $P$ belong to the link of $O$ such as depicted in figure~\ref{fig:linkO}.
\begin{figure}
\centering
\includegraphics[scale=0.25]{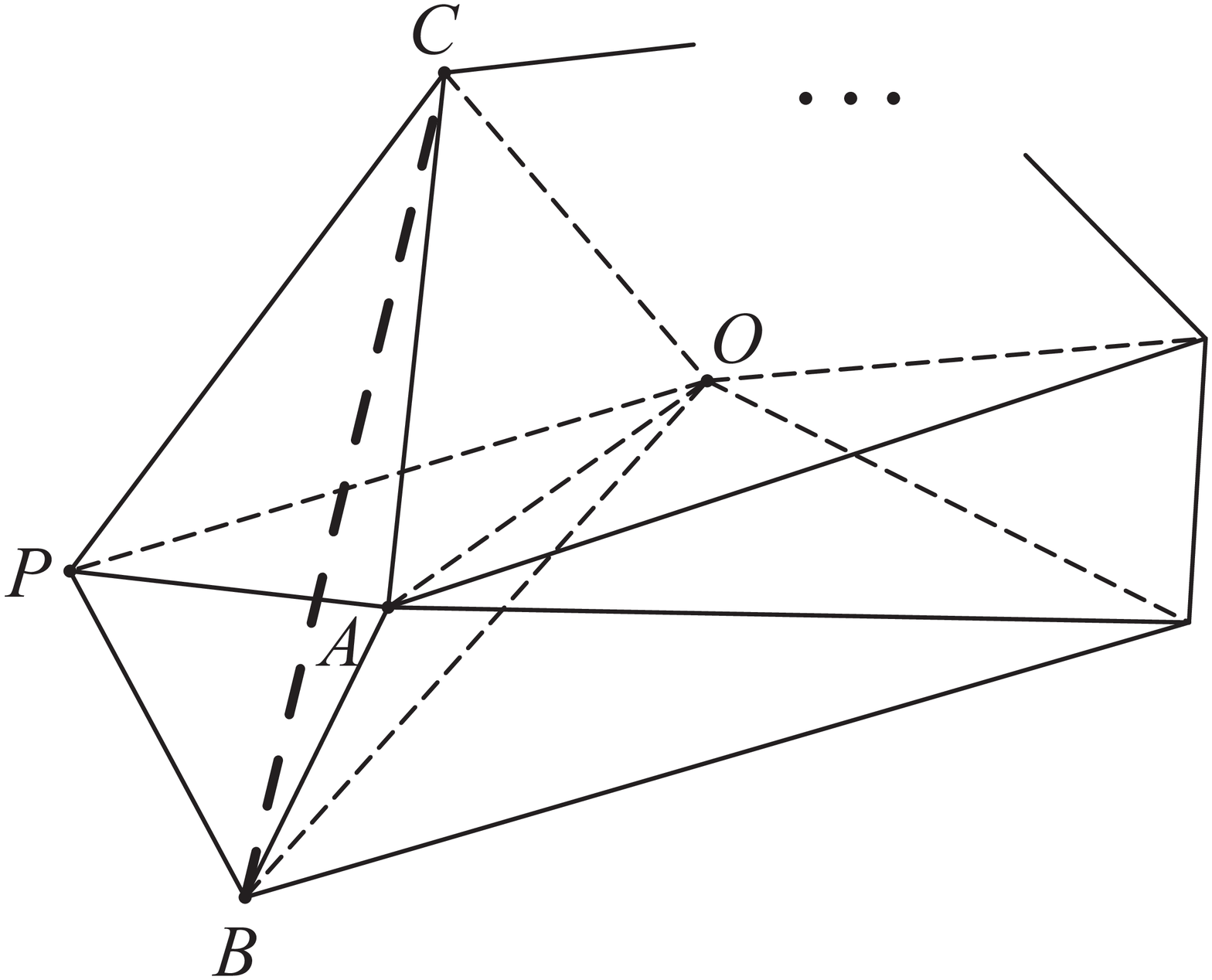}
\caption{The link of vertex $O$}\label{fig:linkO}
\end{figure}
We do a move $2 \to 3$ in the cluster of two tetrahedra $AOPC$ and $BAOP$ adding a new edge $BC$. Doing so several times for other clusters from
$\St(O)$, we can decrease to three the number of edges starting from $A$ and belonging to $\lk(O)$. After that, doing the move $3 \to 2$, we remove
the edge $OA$ from $\St(O)$ and the vertex $A$ from $\lk(O)$. Continuing this way, we decrease to four the number of vertices in $\lk(O)$.

Now the link of $O$ is a tetrahedron and we do the move $4 \to 1$ to remove the vertex $O$ from the triangulation. After that, doing the move $1 \to
4$, we return it back, but in another point of $\mbb{R}^3$, which is required to be in general position with other points (see
remark~\ref{rmk:genpos}). Finally, inverting the whole sequence of the Pachner moves, we return to the initial simplicial complex. It remains only to
note that, by lemma~\ref{lem:pach}, $I_{\rho}(M)$ does not change its value at every step.
\end{proof}

Theorem~\ref{thm:inv} is proven. \qed \bigskip

\section{Calculation of the invariant for $L(p,q)$}\label{sec:calc}

\subsection{Generalities on lens spaces and their triangulations}\label{ssec:lens}
Let $p > q > 0$ be two coprime integers. We identify $S^3$ with the subset $\{(z_1, z_2) \in \mbb{C}^2 \; | \; |z_1|^2 + |z_2|^2 = 1\}$ of
$\mbb{C}^2$. The lens space $L(p,q)$ is defined as the quotient manifold $S^3 / \sim$, where $\sim$ denotes the action of the cyclic group
$\mbb{Z}_p$ on $S^3$ given by:
\begin{equation}\label{18a}
\zeta \cdot (z_1, z_2) = (\zeta z_1, \zeta^q z_2), \quad \zeta = e^{2\pi i/p}.
\end{equation}

As a consequence the universal cover of a lens space is the three-dimensional sphere $S^3$ and
\begin{equation}
\pi_1 \bigl( L(p,q) \bigr) \simeq H_1 \bigl( L(p,q) \bigr) \simeq \mathbb Z_p. \label{eq_Z_p}
\end{equation}

\begin{rmk}
In formula~(\ref{18a}), the group $\mathbb Z_p$ is understood as a multiplicative group of roots of unity of degree~$p$. Below, it will be more
convenient for us to consider it as an additive group consisting of integers between 0 and $p-1$ whose addition is understood modulo~$p$.
\end{rmk}

The full classification of lens spaces is due to Reidemeister and is given in the following theorem.
\begin{thm}[\cite{Reid}]
Lens spaces $L(p, q)$ and $L(p', q')$ are homeomorphic if and only if $p' = p$ and $q' = \pm q^{\pm 1} \bmod p$.
\end{thm}

\begin{figure}
\centering
\includegraphics[scale=0.3]{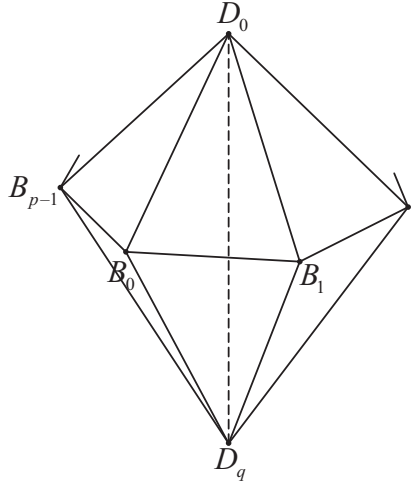}
\caption{A triangulation of $L(p, q)$}\label{fig:complex2}
\end{figure}

Now we describe a triangulation of $L(p,q)$ which will be used in our computations. Consider $p$-gonal bipyramid, i.e. a union of two cones over a
regular polygon of $p$ sides. Denote by $B_0$, $B_1$, \ldots , $B_{p - 1}$ the vertices of the polygon and by $D_0$, $D_q$ the cone apexes. For each
$i$, we glue together the faces $B_i D_0 B_{i + 1}$ and $B_{i + q} D_q B_{i + q + 1}$ (all indices are taken modulo $p$). It is easy to see that the
result is homeomorphic to the lens space $L(p,q)$ (fig.~\ref{fig:complex2}).

Let $g$ be a generator of $\mbb{Z}_p$. Given a Cartesian coordinate system $(x, y, z)$, we define a representation
\begin{equation}
\rho_k \colon g \mapsto
\begin{pmatrix}
\cos \frac{2\pi k}{p} & \sin \frac{2\pi k}{p} & 0 \\
- \sin \frac{2\pi k}{p} & \cos \frac{2\pi k}{p} & 0 \\
0 & 0 & 1
\end{pmatrix}, \quad k = 1, \ldots, p-1,
\label{eq:rhoLpq}
\end{equation}
which sends $g$ in a rotation around the axis $z$ in $\mbb{R}^3$ through the angle $\frac{2\pi k}{p}$, where $k \neq 0$.

We are going to calculate $I_k(L(p, q))$ for the lens space $L(p, q)$ and the representation $\rho_k$ defined by~\eqref{eq:rhoLpq}.

Let us put the vertices of the fundamental family to points $B_0(1, 0, 0)$ and $D_0(\cos \alpha,\allowbreak \sin \alpha, 1)$, where $\alpha =
\frac{\pi}{2} + \frac{\pi k (q-1)}{p}$. Then, one can find the volumes (multiplied by 6) of all consistently oriented tetrahedra, entering in the
triangulation:
\begin{equation}\label{eq:Vi}
\begin{split}
V_i \equiv V_{B_i B_{i+1} D_0 D_q} =
\begin{vmatrix}
x_{D_0} - x_{D_q} & y_{D_0} - y_{D_q} & z_{D_0} - z_{D_q} \\
x_{D_0} - x_{B_i} & y_{D_0} - y_{B_i} & z_{D_0} - z_{B_i} \\
x_{D_0} - x_{B_{i+1}} & y_{D_0} - y_{B_{i+1}} & z_{D_0} - z_{B_{i+1}}
\end{vmatrix} = \\
\begin{vmatrix}
x_{D_0} - x_{D_q} & y_{D_0} - y_{D_q} & 0 \\
x_{B_i} & y_{B_i} & - z_{D_0} \\
x_{B_{i+1}} & y_{B_{i+1}} & - z_{D_0}
\end{vmatrix} = \\
z_{D_0}(x_{D_0} - x_{D_q})(y_{B_{i+1}} - y_{B_i}) - z_{D_0}(y_{D_0}
- y_{D_q})(x_{B_{i+1}} - x_{B_i}) = \\
4 \sin \frac{\pi k}{p} \sin \frac{\pi q k}{p} \cos \frac{2 \pi k i}{p},
\end{split}
\end{equation}
where $i$ runs from 0 to $p-1$ and $B_p = B_0$.

By proposition~\ref{prp:e3rho}, the space $\mf{e}(3)_\rho$ is generated by the differential of a screw motion. Let $d\varphi_z$ and $dz$ be an
infinitesimal rotation around the axis $z$ and infinitesimal translation along this axis respectively. Then, $\mc{C}_0 = \mc{B}_0 = \{d\varphi_z,
dz\}$.

The basis vectors of the space $(dx)$ are the differentials of Euclidean coordinates of $B_0$ and $D_0$. Choose
\[
\mc{B}_1 = \{dz_{B_0}, dx_{D_0}\}.
\]
Then, by direct calculation, we obtain
\begin{equation}
\label{eq:detf1} \det {}_{\mc{B}_1} f_1 = \begin{vmatrix} 0 & 1 \\ y_{D_0} & 0
\end{vmatrix} = - \cos \frac{\pi k (q-1)}{p}.
\end{equation}

The space $(dL)$ has dimension $p+2$ and the set of its basis vectors looks like
\[
\mc{C}_2 = \{dL_{B_0 B_1}, dL_{D_0 B_0}, \ldots, dL_{D_0 B_{p-1}}, dL_{D_0 D_q}\},
\]
Choose
\[
\mc{B}_2 = \{dL_{B_0 B_1}, dL_{D_0 B_{q-1}}, dL_{D_0 B_{p-1}}, dL_{D_0 D_q}\}.
\]
Then,
\begin{equation}\label{eq:detf2}
\begin{split}
\det {}_{\mc{B}_2} f_2 =
\begin{vmatrix}
\pa{L_{B_0 B_1}}{x_{B_0}} & 0 & 0 & 0 \\
\pa{L_{D_0 B_{q-1}}}{x_{B_0}} & \pa{L_{D_0 B_{q-1}}}{y_{B_0}} &
\pa{L_{D_0 B_{q-1}}}{y_{D_0}} & \pa{L_{D_0 B_{q-1}}}{z_{D_0}} \\
\pa{L_{D_0 B_{p-1}}}{x_{B_0}} & \pa{L_{D_0 B_{p-1}}}{y_{B_0}} &
\pa{L_{D_0 B_{p-1}}}{y_{D_0}} & \pa{L_{D_0 B_{p-1}}}{z_{D_0}} \\
0 & 0 & \pa{L_{D_0 D_q}}{y_{D_0}} & 0
\end{vmatrix} = \\
- \pa{L_{B_0 B_1}}{x_{B_0}} \pa{L_{D_0 D_q}}{y_{D_0}}
\begin{vmatrix}
\pa{L_{D_0 B_{q-1}}}{y_{B_0}} & \pa{L_{D_0 B_{q-1}}}{z_{D_0}} \\
\pa{L_{D_0 B_{p-1}}}{y_{B_0}} & \pa{L_{D_0 B_{p-1}}}{z_{D_0}}
\end{vmatrix} =
- \pa{L_{B_0 B_1}}{x_{B_0}} \pa{L_{D_0 D_q}}{y_{D_0}} \pa{(L_{D_0 B_{p-1}} - L_{D_0 B_{q-1}})}{y_{B_0}} z_{D_0} = \\
32 x_{B_0} y_{D_0} z_{D_0} \sin^2 \frac{\pi k}{p} \sin^3 \frac{\pi
k q}{p} (x_{D_0} \cos \frac{\pi k (q-2)}{p} - y_{D_0} \sin \frac{\pi k (q-2)}{p}) = \\
- 32 \sin^2 \frac{\pi k}{p} \sin^3 \frac{\pi k q}{p} \cos \frac{\pi k (q-1)}{p} \sin \frac{\pi k (2q-3)}{p}.
\end{split}
\end{equation}

\begin{thm}\label{thm:ILpq}
The invariant $I_k(L(p, q))$ for the lens space $L(p, q)$ and representation $\rho_k$ defined in~\eqref{eq:rhoLpq} looks as follows
\begin{equation}\label{eq:ILpq}
I_k(L(p, q)) = - \frac{1}{p^2} \left(4 \sin \frac{\pi k}{p} \sin \frac{\pi k q}{p}\right)^4.
\end{equation}
\end{thm}

\begin{proof}
Recall that $\mc{B}_3 = \mc{C}_2 \setminus \mc{B}_2$. In order to find $\det {}_{\mc{B}_3} f_3$, we define a submatrix $f'_3 = {}_{\mc{C}'_2} f_3$ as
a restriction of matrix $f_3$ to the subset
\[
\mc{C}'_2 = \{dL_{D_0 B_i} \mid i=0, \ldots, p-1\}.
\]

Let us write out all the nonzero elements of $f'_3$. Since ${f'_3}^T = f'_3$, we have to write $(f'_3)_{i, j}$ for $j \geq i$ only.

We begin from the diagonal elements. Note that four tetrahedra $B_{i-1} B_i D_0 D_q$, $B_i B_{i+1} D_0 D_q$, $B_{i+q-1} B_{i+q} D_0 D_q$ and $B_{i+q}
B_{i+q+1} D_0 D_q$ share the edge $D_0 B_i$. So, one can use formula~\eqref{eq:4tet} to find $(f'_3)_{i, i} = \frac{1}{l_{D_0 B_i}^2} \pa{\omega_{D_0
B_i}}{l_{D_0 B_i}}$:
\[
(f'_3)_{i, i} = - \left(\frac{V_{D_q B_{i-1} B_i B_{i+1}} V_{B_{i+1} B_{i-1} D_0 D_q}}{V_{i-1} V_i V_{B_{i-1} B_i B_{i+1} D_0}} + \frac{V_{D_q
B_{i+q-1} B_{i+q} B_{i+q+1}} V_{B_{i+q+1} B_{i+q-1} D_0 D_q}}{V_{i+q-1} V_{i+q} V_{B_{i+q-1} B_{i+q} B_{i+q+1} D_0}}\right).
\]
We can simplify this expression using the following obvious equalities:
\[
V_{D_q B_{i-1} B_i B_{i+1}} = V_{B_{i-1} B_i B_{i+1} D_0}, \quad V_{D_q B_{i-1} B_i B_{i+1}} = V_{B_{i-1} B_i B_{i+1} D_0}
\]
and
\[
V_{B_{i+1} B_{i-1} D_0 D_q} = V_{i-1} + V_i, \quad V_{B_{i+q+1} B_{i+q-1} D_0 D_q} = V_{i+q-1} + V_{i+q}.
\]
Then, we obtain
\[
(f'_3)_{i, i} = - \frac{1}{V_{i-1}} - \frac{1}{V_i} - \frac{1}{V_{i+q-1}} - \frac{1}{V_{i+q}}.
\]

Likewise, using~\eqref{eq:skresch}--\eqref{eq:4tet}, one can find the remaining nonzero entries of $f'_3$:
\[
(f'_3)_{i, i+1} = \frac{1}{V_i} + \frac{1}{V_{i+q}},
\]
\[
(f'_3)_{i, i+q-1} = \frac{1}{V_{i+q-1}},
\]
\[
(f'_3)_{i, i+q} = - \frac{1}{V_{i+q-1}} - \frac{1}{V_{i+q}},
\]
\[
(f'_3)_{i, i+q+1} = \frac{1}{V_{i+q}}.
\]

Let us define an auxiliary matrix $\Phi_3 = \Phi_3^T$ with the following nonzero elements
\[
(\Phi_3)_{i, i} = - \frac{1}{V_{i-1}} - \frac{1}{V_i},
\]
\[
(\Phi_3)_{i, i+1} = \frac{1}{V_i},
\]
\[
(\Phi_3)_{i, i-1} = \frac{1}{V_{i-1}}.
\]
Then, the matrix $f'_3$ can be represented as
\begin{equation}
f'_3 = \Phi_3 + E^q \Phi_3 E^{-q} - E^q \Phi_3 - \Phi_3 E^{-q} = (\boldsymbol{1}_p - E^q) \Phi_3 (\boldsymbol{1}_p - E^{-q}),
\end{equation}
where
\begin{equation}\label{eq:matrixE}
E =
\begin{pmatrix}
  0 & 1 & \qquad & \qquad \\
  \qquad & \ddots & \ddots & \qquad \\
  \qquad & \qquad & 0 & 1 \\
  1 & \qquad & \qquad & 0
\end{pmatrix},
\end{equation}
$\boldsymbol{1}_p$ is the identity matrix of size $p \times p$. Units in $E$ form a cyclic diagonal upper the main one.

Also, taking into account that
\[
\Phi_3 = (\boldsymbol{1}_p - E^{-1}) \begin{pmatrix} \frac{1}{V_0} & & \\
& \ddots & \\ & & \frac{1}{V_{p - 1}} \end{pmatrix} (\boldsymbol{1}_p - E),
\]
we finally get
\begin{equation}\label{eq:factorf3}
f'_3 = S_q \, S_{-1} \, R \, S_1 \, S_{-q},
\end{equation}
where we have denoted
\begin{equation}
S_i = \boldsymbol{1}_p - E^i, \qquad R = \diag(V_0^{-1}, V_1^{-1},\ldots V_{p - 1}^{-1}).
\end{equation}

We are going to use the factorization~\eqref{eq:factorf3} of matrix $f'_3$ in order to simplify the matrix ${}_{\mc{B}_3} f_3$ by means of elementary
transformations. The subset $\mc{B}_3$ is chosen in such a way that the matrix $S_q$ in~\eqref{eq:factorf3} loses its $q$th and last rows. Hence, the
$q$th column of this matrix becomes zero. So, we can take away this column and also the $q$th row in $S_{-1}$.

Denote the transformed matrices as $\tilde{S}_q$ and $\tilde{S}_{-1}$. Then, adding the rows in $\tilde{S}_q$, we get it in form
\[
\begin{pmatrix}
1 & {} & {} & -1 \\
{} & \ddots & {} & \vdots \\
{} & {} & 1 & -1
\end{pmatrix}.
\]

Further, let us multiply the matrices $\tilde{S}_q$ and $\tilde{S}_{-1}$ and transform the obtained $(p-2) \times p$ matrix in such a manner that its
first $p-2$ columns form the identity matrix. That is,
\[
\tilde{S}_q \tilde{S}_{-1} \sim \begin{pmatrix} \boldsymbol{1}_{p - 2} & \mb{a} & \mb{b} \end{pmatrix}.
\]
Here, the components of columns $\mb{a}$ and $\mb{b}$ look like
\begin{equation}
\label{eq:ai} \mb{a}_i =
\left\{%
\begin{array}{ll}
i, & i = 1, \ldots, q-1, \\
i - p, & i = q, \ldots, p-2, \\
\end{array}%
\right.
\end{equation}
\begin{equation}
\label{eq:bi} \mb{b}_i = - 1 - \mb{a}_i.
\end{equation}

Now, using the symmetric form~\eqref{eq:factorf3} of the matrix $f'_3$, we find
\begin{equation}
{}_{\mc{B}_3} f_3 \sim \begin{pmatrix} \boldsymbol{1}_{p - 2} & \mb{a} & \mb{b}\end{pmatrix} R \begin{pmatrix} \boldsymbol{1}_{p - 2} \\ \mb{a}^T
\\ \mb{b}^T \end{pmatrix} = R' \begin{pmatrix} \boldsymbol{1}_{p - 2} & \frac{\mb{a}'}{V_{p-2}} & \frac{\mb{b}'}{V_{p-1}}\end{pmatrix} \begin{pmatrix}
\boldsymbol{1}_{p - 2} \\ \mb{a}^T \\ \mb{b}^T \end{pmatrix},
\end{equation}
where $\mb{a}' = R'^{-1} \mb{a}$, $\mb{b}' = R'^{-1} \mb{b}$ and $R' = \diag(V_0^{-1}, V_1^{-1}, \ldots V_{p - 3}^{-1})$.

Denote
\[
P = \begin{pmatrix} \boldsymbol{1}_{p - 2} & \frac{\mb{a}'}{V_{p-2}} & \frac{\mb{b}'}{V_{p-1}}\end{pmatrix} \begin{pmatrix} \boldsymbol{1}_{p - 2} \\
\mb{a}^T \\ \mb{b}^T \end{pmatrix} = \boldsymbol{1}_{p - 2} + \frac{1}{V_{p-2}} \mb{a}' \otimes \mb{a}^T + \frac{1}{V_{p-1}} \mb{b}' \otimes
\mb{b}^T.
\]

It is well known that the eigenvalues of a matrix $\frac{1}{V_{p-2}} \mb{a}' \otimes \mb{a}^T + \frac{1}{V_{p-1}} \mb{b}' \otimes \mb{b}^T$ of rank 2
are
\[
0, \ldots, 0, \lambda_1, \lambda_2,
\]
where $\lambda_1, \lambda_2$ are the eigenvalues of $2\times 2$ matrix:
\[
\begin{pmatrix}
\frac{1}{V_{p-2}} \mb{a}^T \mb{a}' & \frac{1}{V_{p-1}} \mb{a}^T \mb{b}' \\
\frac{1}{V_{p-2}} \mb{b}^T \mb{a}' & \frac{1}{V_{p-1}} \mb{b}^T \mb{b}'
\end{pmatrix}.
\]
This yields
\begin{equation}
\begin{split}
\det P = \frac{1}{V_{p-2} V_{p-1}}
\begin{vmatrix}
V_{p-2} + \mb{a}^T \mb{a}' & \mb{a}^T \mb{b}' \\
\mb{b}^T \mb{a}' & V_{p-1} + \mb{b}^T \mb{b}'
\end{vmatrix} = \\
\frac{1}{V_{p-2} V_{p-1}}
\begin{vmatrix}
V_{p-2} + \mb{a}^T \mb{a}' & V_{p-2} - \sum \mb{a}'_i \\
V_{p-2} - \sum \mb{a}'_i  & 0
\end{vmatrix} =
- \frac{(V_{p-2} - \sum \mb{a}'_i)^2}{V_{p-2} V_{p-1}} = \\
- \frac{1}{V_{p-2} V_{p-1}}(p \sum\limits_{i=q}^{p-1} V_{i-1} - \sum\limits_{i=1}^{p-1} i V_{i-1})^2 = - \frac{4 p^2}{V_{p-2} V_{p-1}} \sin^2
\frac{\pi k q}{p} \sin^2 \frac{\pi k (2q-3)}{p}.
\end{split}
\end{equation}
The second equality holds since $\mb{b}_i = - 1 - \mb{a}_i$. We also took into account that $\sum\limits_{i=0}^{p-1} V_i = 0$, which is true
by~\eqref{eq:Vi}.

So, we find
\[
\det {}_{\mc{B}_3} f_3 = \frac{1}{\prod\limits_{i = 0}^{p-3} V_i} \det P = - \frac{4 p^2}{\prod\limits_{i = 0}^{p-1} V_i} \sin^2 \frac{\pi k q}{p}
\sin^2 \frac{\pi k (2q-3)}{p}.
\]

Combining this with expressions~\eqref{eq:detf1} and~\eqref{eq:detf2} for $\det {}_{\mc{B}_1} f_1$ and $\det {}_{\mc{B}_2} f_2$, we find the torsion
of the complex~\eqref{eq:seq} and then, by formula~\eqref{eq:inv}, the invariant~\eqref{eq:ILpq}. Theorem~\ref{thm:ILpq} is proven.
\end{proof}

\begin{rmk}\label{rmk:reid}
For coprime $p$ and $q$, there are two integers $a$ and $b$ such that $aq + bp = 1$. Let $\zeta$ be a primitive root of unity of degree $p$. Then,
the \textit{Reidemeister torsion} of lens space $L(p, q)$ is defined by the formula (\cite[theorem 10.6]{Tur01}):
\[
\tau_R(L(p, q)) = (1-\zeta^k)^{-1}(1-\zeta^{ka})^{-1},
\]
where $k = 1,\ldots, p-1$. Comparing this with~\eqref{eq:ILpq}, we conclude that for a given $p$ formula~\eqref{eq:ILpq} yields $|\tau_R(L(p,
q))|^{-4}$ up to a constant.
\end{rmk}

\section{A modification of the invariant for lens spaces}\label{sec:modif}

\begin{figure}
\centering
\includegraphics[scale=0.3]{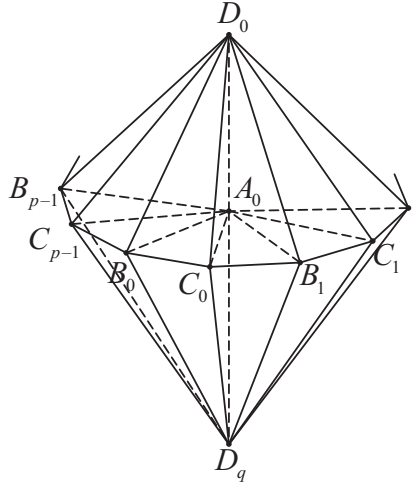}
\caption{Another triangulation of $L(p, q)$} \label{fig:complex}
\end{figure}

Let $M$ be a lens space $L(p,q)$ (see subsection~\ref{ssec:lens}). In figure~\ref{fig:complex}
we depict another triangulation $T$ of $L(p, q)$ which is a subdivision of the one depicted in figure~\ref{fig:complex2}. Now we present $L(p, q)$ as
a triangulated regular $2p$-gonal bipyramid. Take $p - 1$ more copies of this bipyramid (changing indices cyclically modulo $p$) and glue them
together by pairs of faces with the same vertices. Thus we obtain a triangulation of the universal cover $\widetilde{L(p, q)} \cong S^3$. We denote
this triangulation by $\tilde{T}$. It consists of $4p$ vertices, $p(4p+4)$ edges, $8p^2$ 2-simplices and $4p^2$ tetrahedra.

Consider the acyclic complex~\eqref{eq:seq} for $\Gamma(\tilde{T})$:
\begin{equation}\label{eq:acycl1}
0 \xr{} \mf{e}(3) \xr{f_1} (dx) \xr{f_2} (dL) \xr{f_3} (d\Omega) \xr{-f_2^T} (dx)^* \xr{f_1^T} \mf{e}(3)^* \xr{} 0.
\end{equation}
Since $\pi_1(L(p, q)) \simeq \mbb{Z}_p$, all the spaces in~\eqref{eq:acycl1} are finite-dimensional.

Define the following unitary matrix:
\begin{equation}\label{eq:Z}
Z = \frac{1}{\sqrt p}
\begin{pmatrix}
1 & 1 & \cdots & 1 \\
1 & \zeta & \cdots & \zeta^{p - 1} \\
\vdots & \vdots & {} & \vdots \\
1 & \zeta^{p - 1} & \cdots & \zeta^{(p - 1)^2}
\end{pmatrix},
\end{equation}
where $\zeta$ is a primitive root of unity of degree $p$.

In the sequence~\eqref{eq:acycl1}, let us change the bases of $(dx)$ and $(dx)^*$ by means of $Z \otimes \boldsymbol{1}_{12}$ and the bases of $(dL)$
and $(d\Omega)$ by means of $Z \otimes \boldsymbol{1}_{4p+4}$. Then, it is easy to see that the complex~\eqref{eq:acycl1} decomposes into a direct
sum of acyclic complexes:
\begin{equation}\label{eq:acycl2}
\begin{split}
& 0 \xr{} \mf{e}(3) \xr{f_1^{(0)}} (dx)_0 \xr{f_2^{(0)}} (dL)_0 \xr{f_3^{(0)}} (d\Omega)_0 \xr{-(f_2^{(0)})^\dag} (dx)_0^* \xr{(f_1^{(0)})^\dag}
\mf{e}(3) \xr{} 0,
\\
& 0 \xr{} (dx)_j \xr{f_2^{(j)}} (dL)_j \xr{f_3^{(j)}} (d\Omega)_j \xr{-(f_2^{(j)})^\dag} (dx)_j^* \xr{} 0, \qquad j = 1, \ldots, p - 1,
\end{split}
\end{equation}
where $\dag$ means the Hermitian conjugation.

Any matrix element of $f_3^{(j)}$ is, in general, a polynomial in $\zeta^j$ with partial derivatives $\frac{1}{l_a l_b}\frac{\partial
\omega_a}{\partial l_b}$ as its coefficients. These partial derivatives correspond to the five cases considered above, see
formulas~\eqref{eq:skresch}--\eqref{eq:4tet}. However, in our case the matrix $f_3^{(j)}$ has many zero entries. This can be explained in one of two
ways: either the corresponding derivative $\partial \omega_a/ \partial l_b$ vanishes because the edges $a$ and~$b$ do not belong to the same
tetrahedron, or if two edges $a$ and~$b$ belong to the same two-dimensional face (this includes, in particular, the case $a=b$), then the summands in
the derivative
\[
\pa{\omega_a}{l_b} = -\sum_i \frac{\partial (\varphi_a)_i}{\partial l_b},
\]
where $i$ numbers the tetrahedra around $a$, can be grouped in pairs for which the two derivatives $\partial(\varphi_a)_i/\partial l_b$ are equal in
absolute value but differ in signs, because the two corresponding tetrahedra have opposite orientations.

By formula~\eqref{eq:torsion}, the torsions of the complexes~\eqref{eq:acycl2} look as follows
\begin{equation}\label{eq:torsions}
\tau_j = |\det {}_{\mc{B}_1^{(j)}} f_1^{(j)}|^{-2} \, |\det {}_{\mc{B}_2^{(j)}} f_2^{(j)}|^2 \, (\det {}_{\mc{B}_3^{(j)}} f_3^{(j)})^{-1}, \qquad j =
0, \ldots, p - 1,
\end{equation}
where, of course, $f_1^{(j)} = 0$ and $\det {}_{\mc{B}_1^{(j)}} f_1^{(j)} = 1$ for $j \neq 0$.

For simplicity, let us put the vertices $A_i$, $B_i$, $C_i$ and $D_i$ ($i = 0, \ldots, p-1$) to points $A(0, 0, 0)$, $B(1, 0, 0)$, $C(0, 1, 0)$ and
$D(0, 0, 1)$ respectively. Then $V_{ABCD} = 1$ and, by theorem~\ref{thm:inv}, the corresponding invariants look like
\begin{equation}\label{eq:invs2}
 I_j(L(p,q)) = \tau_j
\end{equation}
for $j = 0, \ldots, p - 1$.

Let us introduce the following ordering on the set of $4p+4$ edges in the triangulation $T$ (cf. figure~\ref{fig:complex}):
\begin{equation}\label{eq:ordering}
\begin{split}
& A_0B_0,\ldots, A_0B_{p-1},\; C_0D_0, \ldots, C_{p-1}D_0, \\
& A_0C_0,\ldots, A_0C_{p-1}, \; B_0D_0,\ldots, B_{p-1}D_0, \\
& A_0D_0, A_0D_q,\; B_0C_0, B_1C_0.
\end{split}
\end{equation}
The, we choose
\[
\mc{B}_2^{(0)} = \{1, p+1, 2p+1, 3p+1, 4p+1, 4p+3\}
\]
and
\[
\mc{B}_2^{(j)} = \{1, 2, p+1, p+2, 2p+1, 2p+2, 3p+1, 3p+2, 4p+1, \ldots, 4p+4\}
\]
for $j = 1, \ldots, p-1$. Here, for example, $4p+1$ means $dL_{A_0 D_0}$ and so on.

Due to the acyclicity of~\eqref{eq:acycl2}, we find
\[
\rank f_2^{(0)} = 6, \qquad \rank f_3^{(0)} = 4p - 2
\]
and
\[
\rank f_2^{(j)} = 12, \qquad \rank f_3^{(j)} = 4p - 8
\]
for all $j = 1,\ldots, p - 1$. We see that $\rank f_2^{(j)}$ does not depend on $p$ for all $j$. Therefore, the determinants $\det
{}_{\mc{B}_2^{(j)}} f_2^{(j)}$ in~\eqref{eq:torsions} can be easily found by a direct calculation:
\begin{equation}
\det {}_{\mc{B}_2^{(j)}} f_2^{(j)} =
\begin{cases}
-\frac{1}{p^3}\det {}_{\mc{B}_1^{(0)}} f_1^{(0)}, & j = 0, \\[.4cm]
(\zeta^j  - 1)^5 \, (\zeta^{q j} - 1), & j = 1,\ldots , p - 1,
\end{cases}
\end{equation}
where $\mc{B}_1^{(0)} = \mc{C}_1$.

\begin{thm}\label{thm:ILpq2}
The invariants $I_j(L(p, q))$ look as follows
\begin{equation}\label{eq:ILpq2}
\begin{split}
 I_0(L(p, q)) &= - \frac{1}{p^{12}}, \\
 I_j(L(p,q)) &= \left(4 \sin \frac{\pi j}{p} \sin \frac{\pi q j}{p}\right)^6.
\end{split}
\end{equation}
\end{thm}

\begin{proof}
Recall that ${\mc B}_3^{(j)} = \mc{C}_3\setminus\mc{B}_2^{(j)}$. Let us find $\det {}_{\mc{B}_3^{(j)}} f_3^{(j)}$. The ordering~\eqref{eq:ordering}
gives the following block structure for matrix $f_3^{(j)}$ (here and below the empty spaces are of course occupied by zeroes):
\begin{equation}\label{eq:f3j}
f_3^{(j)} =
\begin{pmatrix}
   \begin{matrix}
   \boldsymbol{0}_p & S_j  \\
   S_j^{\dag} & \boldsymbol{0}_p
   \end{matrix} & \vline &  {} & \vline \\
\hline
   {} & \vline & \begin{matrix}
   \boldsymbol{0}_p & R_j  \\
   R_j^{\dag} & \boldsymbol{0}_p
\end{matrix} & \vline \\
\hline
   {} & \vline &  {} & \vline & \begin{matrix}
   0 & 0 & \sigma_j  &  - \sigma_j \\
   {0} & {0} & { - \sigma_j } & \sigma_j \\
   \sigma_j  &  - \sigma_j  & 0 & 0 \\
    - \sigma_j  & \sigma_j  & 0 & 0
 \end{matrix}
\end{pmatrix},
\end{equation}
where
\[
\sigma_j = \sum\limits_{i = 0}^{p - 1}{\zeta^{ij}} =
\begin{cases}
 p, & j = 0, \\
 0, & j \neq 0.
\end{cases}
\]
The matrices $S_j$ and $R_j$ have size $p\times p$ and can be represented as follows:
\begin{equation}\label{eq:SjRj}
\begin{split}
S_j &= (\boldsymbol{1}_p - \zeta^{qj} E^{-q}) \, (\boldsymbol{1}_p - E^{-1}), \\
R_j &= (\boldsymbol{1}_p - \zeta^{qj} E^{-q}) \, (\boldsymbol{1}_p - E),
\end{split}
\end{equation}
where the matrix $E$ is defined in~\eqref{eq:matrixE}. It follows from the acyclicity and~\eqref{eq:f3j} that the rank of $S_j$ and $R_j$ equals $p -
1$ for $j = 0$ and $p - 2$ for $j = 1,\ldots , p - 1$.

\begin{lem}
Let $\mc{D}_j$ and $\mc{E}_j$ be the corresponding subsets of $\mc{B}_3^{(j)}$ for the matrices $S_j$ and $R_j$ respectively. Then,
\begin{equation}\label{eq:detSj}
\det {}_{\mc{D}_j} S_j =
\begin{cases}
p, & j = 0, \\[.4cm]
\frac{\zeta^j - 1}{\zeta^{q j} - 1} , & j = 1, \ldots, p - 1
\end{cases}
\end{equation}
and
\begin{equation}\label{eq:detRj}
\det {}_{\mc{E}_j} R_j =
\begin{cases}
p, & j = 0, \\[.4cm]
\frac{\zeta^{-j} - 1}{\zeta^{q j} - 1} , & j = 1, \ldots, p - 1.
\end{cases}
\end{equation}
\end{lem}

\noindent \textit{Proof} is done by similar methods as we used proving theorem~\ref{thm:ILpq}. \qed \bigskip

Now, according to~\eqref{eq:f3j},
\begin{equation}
\det {}_{\mc{B}_3^{(j)}} f_3^{(j)} =
\begin{cases}
 - |\det {}_{\mc{D}_0} S_0|^2 \, |\det {}_{\mc{E}_0} R_0|^2 \, p^2 = - p^6, & j = 0, \\[.4cm]
 |\det {}_{\mc{D}_j} S_j|^2 \, |\det {}_{\mc{E}_j} R_j|^2 = \frac{\sin^4 \frac{\pi j}{p}}{\sin^4 \frac{\pi q j}{p}}, & j = 1,\ldots , p - 1.
\end{cases}
\end{equation}
Finally, we calculate the torsions~\eqref{eq:torsions} which are equal to the invariants, according to~\eqref{eq:invs2}. Thus we
get~\eqref{eq:ILpq2}.
\end{proof}

\begin{rmk}
Since we performed the unitary change-of-basis transformations, then the torsion of~\eqref{eq:acycl1} does not change. It follows that multiplying
together the found values~\eqref{eq:ILpq2}, we get the invariant for $\widetilde{L(p ,q)} \cong S^3$, that is
\[
\prod\limits_{j=0}^{p-1} I_j(L(p, q)) = I(S^3) = -1.
\]
\end{rmk}

\bibliographystyle{amsplain}

\end{document}